\documentclass[conference]{IEEEtran}
\usepackage{cite}
\ifCLASSINFOpdf
\else
\fi
\usepackage[cmex10]{amsmath}
\usepackage{amsfonts}
\newtheorem{defi}{Definition}

\newtheorem{thm}[defi]{Theorem}
\newtheorem{lem}[defi]{Lemma}

\newtheorem{prop}[defi]{Proposal}

\begin{document}
\title{Mutually unbiased bases as submodules and subspaces}

\author{\IEEEauthorblockN{Joanne L. Hall\IEEEauthorrefmark{1} and Jan \v{S}\v{t}ov\'\i\v{c}ek\IEEEauthorrefmark{2}}
\IEEEauthorblockA{Department of Algebra\\
Charles University in Prague\\
186 75 Praha 8, Sokolovska 83, Czech Republic \\
\IEEEauthorrefmark{1}Email: hall@karlin.mff.cuni.cz\\
\IEEEauthorrefmark{2}Email: stovicek@karlin.mff.cuni.cz}
%\and
%\IEEEauthorblockN{}
%\IEEEauthorblockA{Department of Algebra\\
%Charles University in Prague\\
%186 75 Praha 8, Sokolovska 83, Czech Republic \\
%}
}

%Telephone: (800) 555--1212, Fax: (888) 555--1212}


\maketitle

\begin{abstract}
%\boldmath
  Mutually unbiased bases (MUBs) have been used in  several cryptographic and communications applications.  
There has been much speculation regarding connections between MUBs and  finite geometries.  Most of which has focused on a connection with projective and affine planes.  
We propose a connection with higher dimensional projective geometries and projective Hjelmslev geometries.  We show that this proposed geometric structure is present in several constructions of MUBs.
\end{abstract}


% \fi
\IEEEpeerreviewmaketitle

\section{Introduction}
Mutually unbiased bases (MUBs) are  a structure first defined in a quantum physics context in 1960 \cite{Schw60}.  Since then MUBs have been used in  quantum key distribution protocols \cite{BB84, SARG04}, and can be used to construct signal sets for communications systems \cite{Alltop80, DY07}.  

A basis for $\mathbb{C}^d$ is \emph{orthonormal} if all basis vectors are orthogonal and of unit length.   Two orthonormal bases $\mathcal{B}_0$ and $\mathcal{B}_1$ in $\mathbb{C}^d$ are called \emph{mutually unbiased} if   $|\langle
\phi|\psi\rangle|^2=1/d$ for all $\phi \in \mathcal{B}_0$ and
$\psi\in \mathcal{B}_1$.

The maximum number of mutually unbiased bases in $\mathbb{C}^d$ is $d+1$ \cite{WF89}.  A set of $d+1$ MUBs is called \emph{complete}, it is complete sets of MUBs that are of most use in the communications applications.  While constructions of complete sets of MUBs in $\mathbb{C}^d$ are known when $d$ is a prime power \cite{WF89}, it is unknown if such  complete sets exist in non-prime power dimensions. 

There has been much speculation regarding connections between MUBs and  finite geometries \cite{SPR04, Beng05,  SP06, Woot06}.  Most of this has focused on a connection with projective and affine planes.  

%\begin{conj}[SPR]\cite{SPR04} A complete set of MUBs exists in $\mathbb{C}^d$ if and only if a projective plane of order $d$ also exits.
%\end{conj}
%Note that a projective plane of order $d$ exists if and only if an affine plane of order $d$ exists.

The evidence for connections between MUBs and finite geometries falls into two categories:  counting arguments \cite{SPR04, SP06}, and structures which construct both MUBs and finite geometries.  These structures include planar functions \cite{RS07,KR03},  symplectic spreads \cite{Kantor2011} as well as specific affine planes \cite{GHW04, RDH10}.  
%Evidence against a strong connection uses structures which are similar to affine planes, but have no corresponding structures similar to MUBs \cite{WD2010}.   
%\marginpar{change}

We investigate higher dimensional projective geometries and show that some sets of MUBs may be regarded as subspaces.  Note that in order for these higher order projective geometries to exist, a projective  plane of the appropriate size must also exist.  If all MUBs are subspaces of larger projective geometries, then a connection between MUBs and projective planes would be proven.  Alas we do not go so far.

%\marginpar{change} 
It has been shown that complete sets of MUBs are equivalent to orthogonal decompositions of the Lie algebra  $sl_n(\mathbb{C})$ \cite{BSTW07}, however finding orthogonal decompositions of Lie algebras is as difficult a task as finding sets of MUBs. Some work has been done classifying  Lie Algebras using  projective geometry \cite{LM2002}, but these results have as yet not been applied to decompositions of $sl_n(\mathbb{C})$. 

Some sets of MUBs have been show to have an Abelian group structure \cite{KMR06, Hall11}.  We go further by showing that some complete sets of MUBs may be regarded as submodules of the appropriate free module, and   as subspaces of a projective geometry over that module. %(for odd dimensions) and subspaces of projective Hjelmslev geometry (for even dimensions).  
%\marginpar{change}

\section{Preliminaries}
\subsection{Constructions of MUBs \label{sec:constructions}}
We investigate three non-equivalent constructions of MUBs. This first construction is based on planar functions over a finite field.  For more on planar functions see for example \cite{CM97D}. Let $\omega_p=e^{\frac{2i\pi}{p}}$.
\begin{thm}[Planar function construction]\label{thm:RS07}\cite[Thm 4.1]{RS07}
Let $\mathbb{F}_q$ be a field of odd characteristic $p$.  Let $\Pi(x)$ be a planar function on $\mathbb{F}_q$.  Let $V_a=\{v_{ab}:b\in \mathbb{F}_q\}$ be the set of vectors
\begin{equation}\label{eqn:planarv}
\vec{v}_{ab}= \frac{1}{\sqrt{q}}\left(\omega_p^{\mbox{\emph{tr}}(a\Pi(x)+bx)}\right)_{x\in\mathbb{F}_q}\
\end{equation}
with $a,b\in\mathbb{F}_q$. The standard basis $E$ along with the sets $V_a$, $a\in \mathbb{F}_q$, form a complete set of $q+1$ MUBs  in $\mathbb{C}^q$.
\end{thm}
The following construction has been shown to be equivalent to the planar function construction when using $\Pi(x)=x^2$ \cite{GR09}.  We highlight it as  the submodule and subspaces structure appear in a  different way to the planar function construction.
\begin{thm}[Alltop Construction]\cite{Alltop80}\cite[Thm 1]{KR03}\label{alltop}
Let $\mathbb{F}_q$ be a finite field of odd characteristic $p\geq5$.  Let $V_a=\{\vec{v}_{ab}:b\in \mathbb{F}_q\}$ be the set of vectors
\begin{equation}\label{Alltopv}
\vec{v}_{ab}= \frac{1}{\sqrt{q}}\left(\omega_p^{\mbox{\emph{tr}}((x+a)^3+b(x+a))}\right)_{x\in\mathbb{F}_q}
\end{equation}
with $a,b\in\mathbb{F}_q$.   The standard basis $E$ along with the sets $V_a$, $a\in \mathbb{F}_q$, form a complete set of $q+1$ MUBs  in $\mathbb{C}^q$.
\end{thm}
The next construction stems from a symplectic spread.
\begin{thm}\cite[3.5(b)]{Kantor2011}\label{thm:sym}
Let $\mathbb{F}_{p^n}$ be a field of odd characteristic $p$, with $n$ odd.  Let $s$ and $n$ be coprime, such that $s< n/2$. Let $V_a=\{v_{ab}:b\in \mathbb{F}_q\}$ be the set of vectors
\begin{equation}\label{eqn:symv}
\vec{v}_{ab}= \frac{1}{\sqrt{q}}\left(\omega_p^{\mbox{\emph{tr}}(ax+bx^{p^{n-s}+1}+b^{p^s}x^{p^s+1})}\right)_{x\in\mathbb{F}_q} 
\end{equation}
with $a,b\in\mathbb{F}_q$. The standard basis $E$ along with the sets $V_a$, $a\in \mathbb{F}_q$, form a complete set of $q+1$ MUBs  in $\mathbb{C}^q$.
\end{thm}
The next construction uses Galois rings.
\begin{thm}[Galois ring construction]\label{even}\cite[Thm 3]{KR03}
Let $GR(4,n)$ be  Galois ring of characteristic $4$ and  Teichm\"{u}ller set $\mathcal{T}_n$.  Let $i=\omega_4=\sqrt{-1}$.  Let $V_a=\{\vec{v}_{ab}:b\in \mathcal{T}_r\}$ be the set of vectors
\begin{equation}\label{eqn:ringmubs}
\vec{v}_{ab}=\frac{1}{\sqrt{2^n}}\left(i^{\mbox{\emph{{tr}}}[(a+2b)x]}\right)_{x\in\mathcal{T}_n}
\end{equation}
$a,b\in \mathcal{T}_n$.  The standard basis $E$ along with the sets $V_a$, $a\in\mathcal{T}_n$, form a complete set of $2^n+1$ MUBs in $\mathbb{C}^{2^n}$.
\end{thm}
These are not the only known constructions of complete sets of MUBs \cite{Kantor2011}, but are good starting point for an investigation.

\subsection{Algebraic Structures}
Let $R$ be a ring with unity, a left $R$-module  is an Abelian group, $M$, together with a product $R\times M\mapsto M$ which satisfies the following: for all $r_i,r_2\in R$ and $a_i,a_2\in M$
\begin{align}
1a= & a,\\
(r_1r_2)a= & r_1(r_2 a)\\
(r_1+r_2)a = & r_1a+r_2a\\
r(a_1+a_2) = & ra_1+ra_2
\end{align}
This is familiar as the left  axioms of a vector space. All $\mathbb{F}$-modules where $\mathbb{F}$ is a field are vector spaces.  Theorem \ref{even} uses a ring to construct MUBs, hence we need the more general object of a module.  We are only concerned with commutative rings, thus all modules in consideration are both left and right modules.  An (left and right) $R$ module is \emph{free} if it is isomorphic to $R^d$ for some $d$.  

The trace map, familiar from finite fields, may also be used in Galois rings \cite[\S 14]{Wan2012}. Properties of  trace map for $GR(4,n)$ have been well studied in a coding theory context \cite{NK97}.     
\begin{thm}\cite[Thms 7.12, 14.34, 14.37]{Wan2012}\label{thm:trace}
The trace map,  $\mbox{\emph{tr}}:GR(p^s,n)\mapsto GR(p^s,1)$ has the following properties: 
\begin{enumerate}
\item
For all $r\in GR(p^s,1)$ and $x\in GR(p^s,n)$, $r\mbox{\emph{tr}}(x)=\mbox{\emph{tr}}(rx)$. \label{part:r}
\item
$\mbox{\emph{tr}}(\alpha)=0$ if an only if there exists $\beta\in R'$ such that $\alpha=\beta-\phi(\beta)$.\label{part:0}
\end{enumerate}
where $\phi$ is the generalized Frobenius automorphism.
Note that $GR(p^1,n)\cong \mathbb{F}_{p^n}$.
\end{thm}
For further on Galois rings and fields we refer the reader to \cite{Wan2012}.

\subsection{Geometric Structures}
The geometric structures we are investigating are projective geometries,  $PG(d-1,q)$, defined over a finite field   and projective Hjelmslev geometries $PHG(d-1,GR(4,1))$, defined over a Galois ring.

Let $M$ be an $R$ module that is a submodule of $R^d$.  If $R$ is a field, then any submodule  is a subspace of $R^d$.  If $R$ is a Galois ring then any free submodule is a subspace of $R^d$ \cite{Kreuzer91}.

%If $M$ is a free direct summand of $R^d$ then $M$ is a \emph{subspace} of $R^d$, ie. %$M$ is a free $R$-module, and there exists another $R$-module, $L$ such that %$R^d=M\oplus L$ \cite[\S 1]{Veld95}.

 \begin{defi}\label{defi:PG}
The projective geometry constructed from $\mathbb{F}_q$,  $PG(d-1,q)$ is the set of subspaces of $\mathbb{F}_q^d$.  
$\langle\vec{x}\rangle$ is a point of $PG(d-1,q)$ and  represents all  vectors $\rho\vec{x}$ in $\mathbb{F}_q^d$ such that $\rho\in\mathbb{F}_q^*$ and at least one of the entries of $\vec{x}$ is non-zero. 
 \end{defi}
 
\begin{defi}\label{defi:PHG}\cite{Veld95}
The projective Hjelmslev geometry constructed from $GR(4,1)$,  $PHG(d-1,GR(4,1))$ is the set of subspaces of $GR(4,1)^d$.  
$\langle\vec{x}\rangle$ is a point of $PHG(d-1,GR(4,1))$ and  represents all  vectors $\rho\vec{x}$ in $GR(4,1)^d$ such that $\rho$ is a unit of $GR(4,1)$ and at least one of the entries of $\vec{x}$ is a unit of $GR(4,1)$. 
 \end{defi}
Note that $PG(d-1,q)\cong PHG(d-1,\mathbb{F}_q)$.

\section{MUBs as submodules and subspaces}
\subsection{Conjecture}
\begin{prop}\label{prop}
Let $X$ be a complete set of MUBs which contains the standard basis in $\mathbb{C}^d$.  Let $N$ be the set containing all the vectors from $X$, except the standard basis vectors.  Let the vectors in $N$ be of the form $\alpha\omega_q^{\vec{x}}$ where $\alpha\in\mathbb{R}$, $\omega_q$ is a $q^{th}$ root of unity, and $\vec{x}\in \mathbb{Z}_q^d$.  Let $\odot$ represent component wise multiplication,  let \begin{equation}
\vec{v}\hat{\odot}\vec{u}=\frac{\vec{v}\odot\vec{u}}{|\vec{v}\odot\vec{u}|}
\end{equation}
and let $N'=\{\vec{u}\hat{\odot}\vec{v}^*:\vec{u},\vec{v}\in N\}$,   $M=\{\vec{x}:\alpha\omega_q^{\vec{x}}\in N\}$, and $M'=\{\vec{x}-\vec{y}:\vec{x},\vec{y}\in M\}$.  Let $U'\subset M'$ be the set containing the vectors from $M'$ for which every entry is a non-unit, then
\begin{enumerate}
\item $N'$ is a $\mathbb{Z}_q$-module.
\item $M'\setminus U'$ is the set of vectors representing a subspace of a projective geometry over $\mathbb{Z}_q$.
\end{enumerate}
\end{prop}
We show this proposal is true for each of the constructions of MUBs mentioned in section \ref{sec:constructions}.  This proposal says nothing about the existence of MUBs which are not constructed from a ring.  All projective geometries and projective Hjelmslev geometries of dimension greater than 2 have an algebraic structure \cite[\S 1.4]{dem97},\cite{Kreuzer91}.  It may be the same for complete sets of MUBs.

%\marginpar{change}
MUBs for which the set of vectors forms a group under point-wise multiplication have been studied \cite{GR09}.  Our construction is more general in that the algebraic structure is in the set of vectors generated by point-wise multiplication.

\subsection{Counting}
Much of the evidence for connections between MUBs and geometric structures stems from similarities in cardinality.  We show that Proposal \ref{prop} is plausible in general by using cardinalities. 
\begin{lem}\label{lem:oddcount}
Let $q=p^n$, with $p$ odd, each point in $PG(q-1,p)$ is represented by $p-1$ vectors.  The number of vectors represented by the points in a $(2n-1)$-dimensional subspace of $PG(q-1,q)$, with the addition of $\vec{0}$ is the same as the number of vectors in a complete set of MUBs in $\mathbb{C}^q$ minus the standard basis.
\end{lem}
\begin{IEEEproof}
Let $X$ be an $m$ dimensional subspace of $PG(p^n-1,p)$ then there are $\frac{p^{m+1}-1}{p-1}$ points, each of which may be represented by $p-1$ different vectors.  Add the vector $\vec{0}$. $(p-1)\frac{p^{m+1}-1}{p-1}+1=p^{m+1}$. The number of vectors in a complete set of MUBs in $\mathbb{C}^q$, minus the standard basis is $p^{2n}$.  Thus if we require every vector in the set of MUBs to represent a point in the subspace, we need a $2n-1$ dimensional subspace of  $PG(p^n-1,p)$.  
\end{IEEEproof}

\begin{lem}\label{lem:evencount}
Each point in $PHG(2^n-1, GR(4,1))$ is represented by $2$ vectors.  The number of vectors represented by the points in a $2^{n-1}$ dimensional subspace of $PHG(2^n-1,GR(4,1))$, with the addition of $2^n$ vectors containing no unit elements  is the same as the number of vectors in a complete set of MUBs in $\mathbb{C}^{2^n}$ without the standard basis.
\end{lem}
\begin{IEEEproof}
Let $X$ be an $m$ dimensional subspace of $PHG(2^n-1,GR(4,1))$ then there are $2^m$ points in each of $2^{m+1}-1$ neighbourhoods, each of which may be represented by $2$ different vectors.  $2.2^m(2^{m+1}-1)=2^{2(m+1)}-2^{m+1}$, which, when we add $2^m$ vectors which are generated by non units, is the number of vectors in a complete set of MUBs in $\mathbb{C}^{2^{m+1}}$, minus the standard basis.
\end{IEEEproof}

\subsection{Odd dimensions}
We now show that for specific families of  MUBs proposal \ref{prop} is true.
\begin{thm}\label{thm:odd}
Let $X$ be the complete set of MUBs in $\mathbb{C}^{p^n}$ generated by the planar function construction (Thm \ref{thm:RS07}).  Let $N\subset X$ be the set of vectors $\vec{X}=\frac{1}{\sqrt{d}}\omega_p^{\vec{x}}$ where $\vec{x}\in \mathbb{F}_p^{p^r}$.  Let $M=\{\vec{x}:\omega_p^{\vec{x}}\in N\}$, then
\begin{enumerate}
\item $\langle N,\hat{\odot}\rangle$ is an $\mathbb{F}_p$-module. 
\item  $M$ is a $2n-1$ dimensional subspace  of    $PG(p^{n}-1,p)$.
\end{enumerate}
\end{thm}
\begin{IEEEproof}
1. 
 Let $\vec{v}_{ab}$ and $\vec{v}_{cd}$ be given as in equation (\ref{eqn:planarv}).  
\begin{equation} \vec{v}_{ab}\hat{\odot}\vec{v}_{cd}= \frac{1}{\sqrt{q}}\left(\omega_p^{\mbox{tr}[(a+c)\Pi(x)+(b+d)x]}\right)_{x\in\mathbb{F}_q}
\end{equation}
with $a,b,c,d\in\mathbb{F}_q$. Hence $(\vec{v}_{ab}\hat{\odot}\vec{v}_{cd})\in N$,  $\vec{v}_{00}$ acts as an identity element, with $\vec{v}_{ab}\hat{\odot}\vec{v}_{((-a)(-b)}=\vec{v}_{00}$ ensuring every element has an inverse; commutativity comes from $\mathbb{F}_q$. Thus we have shown that $\langle N,\hat{\odot}\rangle$ is an Abelian group (See also \cite[Lem 2.84]{Hall11}).  To show that it is a module $\mathbb{F}_p\times N\mapsto N$, let $r\in\mathbb{F}_p$, Let $\star$ be an operation on the set $N$ which corresponds to scalar multiplication on the set $M$. 
\begin{eqnarray}
r\star\vec{v}_{ab}= & \frac{1}{\sqrt{q}}\left(\omega_p^{r\mbox{tr}\left(a\Pi(x)+bx\right)}\right)_{x\in\mathbb{F}_q}
\end{eqnarray}
with $a,b\in\mathbb{F}_q$.   By Theorem \ref{thm:trace} 
\begin{eqnarray}
r\star\vec{v}_{ab}= & \frac{1}{\sqrt{q}}\left(\omega_p^{\mbox{tr}\left(ra\Pi(x)+rbx\right)}\right)_{x\in\mathbb{F}_q}
\end{eqnarray}
with $a,b\in\mathbb{F}_q$. Hence for all $r\in\mathbb{F}_p$ and $\vec{v}_{ab}\in N$,  $r\star \vec{v}_{ab}\in N$.  The properties of $\mathbb{F}_p$ ensure that the module axioms are satisfied.
  
2. 
Part 1. shows that $M$ is a submodule, and thus forms a subspace of $\mathbb{F}_{p}^{p^n}$.  The counting results of Lemma 10 show the size of the subspace.
\end{IEEEproof}
For all $a,b,c,d\in\mathbb{F}_q$, any element in $\vec{v}_{ef}\in N$ can be constructed as $\vec{v}_{ef}=\vec{v}_{ab}\hat{\odot}\vec{v}^*_{cd}$ for some $\vec{v}_{ab},\vec{v}_{cd}\in N$.  Thus in the definition of Proposal \ref{prop}, $N=N'$ and  $M=M'$. Hence Proposal \ref{prop} holds for planar function MUBs.

\begin{thm}\label{thm:alltopmodule}
Let $X$ be the complete set of MUBs in $\mathbb{C}^{p^n}$ generated by the Alltop construction (Thm \ref{alltop}). 
Let $S\subset X$ be the set of vectors $\vec{X}=\frac{1}{\sqrt{d}}\omega_p^{\vec{x}}$ where $\vec{x}\in \mathbb{F}_p^{p^r}$.  Let $T=\{\vec{x}:\omega_p^{\vec{x}}\in N\}$, 
 \mbox{$S'=\{\vec{v}\hat{\odot}\vec{u}:\vec{v},\vec{u}\in S\}$} and $T'=\{\vec{x}+\vec{y}:\vec{x},\vec{y}\in T\}$, then
\begin{enumerate}
\item $\langle S',\hat{\odot}\rangle $ is an $\mathbb{F}_p$-module. 
\item  $T'$ is a $2n-1$ dimensional subspace  of    $PG(p^{n}-1,p)$.
\end{enumerate}
\end{thm}
\begin{IEEEproof}
Let $\vec{v}_{ab},\vec{v}_{cd}$ be as defined in equation (\ref{Alltopv}).  
We now show that $S'=N$ and $T'=M$, with $M,N$  from Theorem \ref{thm:odd}.
\begin{align}\label{Alltopinnerproduct} &  \vec{v}_{ab}\hat{\odot}\vec{v}^*_{cd}=\nonumber\\
& \frac{1}{q}\left( \omega_p^{3(a-c)x^2+(3a^2\!-3c^2+b-d)x+(a^3\!-c^3+ba-dc)}\right)_{x\in\mathbb{F}_q}
\end{align}
which is a quadratic in $x$, and hence a planar function.  Theorem \ref{thm:odd} may be used.
\end{IEEEproof}
This  highlights that structures which are  not  present in sets of vectors,  may be present in another way, see also \cite{RDH10}.
We use essentially the same proof for  the MUBs generated by Theorem \ref{thm:sym}.
\begin{thm}\label{thm:sym_module}
Let $X$ be the complete set of MUBs in $\mathbb{C}^{p^n}$ generated by the construction of Theorem \ref{thm:sym}.  Let $N\subset X$ be the set of vectors $\vec{X}=\frac{1}{\sqrt{d}}\omega_p^{\vec{x}}$ where $\vec{x}\in \mathbb{F}_p^{p^r}$.  Let $M=\{\vec{x}:\omega_p^{\vec{x}}\in N\}$ then 
\begin{enumerate}
\item $\langle N,\hat{\odot}\rangle$ is a  $\mathbb{F}_p$-module.
\item $M$ is a $2n-1$ dimensional subspace of $PG(p^n-1,p)$.
\end{enumerate}
\end{thm}
\begin{IEEEproof}
1. We use the  operations $\hat{\odot}$ and $\star$ as in Theorem \ref{thm:odd}. From equation \ref{eqn:symv},
\begin{align}\label{eqn:symdagger}
& \vec{v}_{ab}\hat{\odot} \vec{v}^*_{cd}=\nonumber\\ &  \frac{1}{\sqrt{q}}\left(\omega_p^{\mbox{tr}((a-c)x+(b-d)x^{p^{n-s}+1}+(b^{p^s}-d^{p^s})x^{p^s+1})}\right)_{x\in\mathbb{F}_q}
\end{align}
with $a,b,c,d\in\mathbb{F}_q$.  Let $\phi(b)=b^p$ be the Frobenius automorphism \cite[\S 7.1]{Wan2012}, then 
$b^{p^s}=\phi^{p^{s-1}}(b)$ and hence  $\phi^{p^{s-1}}(b-d)=\phi^{p^{s-1}}(b)-\phi^{p^{s-1}}(d)$.  Using this fact we can rearrange equation (\ref{eqn:symdagger})
\begin{align}\label{eqn:symdagger2}
& \vec{v}_{ab}\hat{\odot} \vec{v}_{cd}^*=\nonumber\\ & \frac{1}{\sqrt{q}}\left(\omega_p^{\mbox{\emph{tr}}((a-c)x+(b-d)x^{p^{n-s}+1}+(b-d)^{p^s}x^{p^s+1})}\right)_{x\in\mathbb{F}_q}
\end{align}
  with $a,b\in\mathbb{F}_q$. Showing that $\vec{v}_{ab}\hat{\odot} \vec{v}_{cd}^*\in N$.  
As with Theorem \ref{thm:odd}, we use the operation $\star$ and see that  $N$ is an $\mathbb{F}_p$ module.

2. The proof is the same as for Theorem \ref{thm:odd}.
\end{IEEEproof}
As with Theorem \ref{thm:odd}, we find that $M=M'$ and $N'=N$ for $M,M',N,N'$  as in Proposal \ref{prop}.

These three structures based on finite fields all conform to the structure of Proposal \ref{prop}.  %Note that the symmplectic spreads and planar functions  form both  affine planes as well as  sets of MUBs \cite{dem97,CCKS97}.  
%\marginpar{change}

\subsection{Even dimensions}

\begin{thm}\label{thm:even}
Let $X$ be the complete set of MUBs in dimension $d=2^n$ generated by the Galois ring   construction \cite{KR03}.  Let $N\subset X$ be the set of vectors $\vec{X}=\frac{1}{\sqrt{d}}i^{\vec{x}}$ where $\vec{x}\in GR(4,1)^{2^n}$.  Let $M=\{\vec{x}:i^{\vec{x}}\in N\}$, then 
\begin{enumerate}
\item $N$ is a  $GR(4,1)$-module.
\item $M$ is a $2^{n-1}$ dimensional subspace of $PHG(2^{n}-1,GR(4,1))$.
\end{enumerate}
\end{thm}
\begin{IEEEproof}
1. 
Let $\alpha=a+2b$, and $\beta=c+2d$ where $a,b,c,d\in\mathcal{T}_n$ the Teichmuller set of $GR(4,n)$.  Then equation (\ref{eqn:ringmubs}) becomes
\begin{equation}\label{eqn:ringmubsalpha}
\vec{v}_{\alpha}=\frac{1}{\sqrt{2^r}}\left(i^{\mbox{\emph{{tr}}}[\alpha x]}\right)_{x\in\mathcal{T}_n}
\end{equation}
$\alpha\in GR(4,n)$.  Let $\hat{\odot}$ be as in Proposal \ref{prop}
\begin{eqnarray}
\vec{v}_\alpha\hat{\odot}\vec{v}_\beta= &  \frac{1}{\sqrt{2^r}}\left(i^{\mbox{\emph{{tr}}}[\alpha+\beta x]}\right)_{x\in\mathcal{T}_r}
\end{eqnarray}
$\vec{v}_\alpha\hat{\odot}\vec{v}_\beta\in M$.  $\vec{v}_{0}$ is the identity, $\vec{v}_{\alpha}\hat{\odot}\vec{v}_{-\alpha}=\vec{v}_{0}$ showing inverses, and commutativity is given by the properties of Galois rings.

Let $\star$ be the operation $GR(4,1)\times N$ that corresponds to scalar multiplication on $M$.
\begin{align}
r\star \vec{v}_\alpha= & \frac{1}{\sqrt{2^n}}\left(i^{r\mbox{\emph{{tr}}}[\alpha x]}\right)_{x\in\mathcal{T}_r}\nonumber\\
= & \frac{1}{\sqrt{2^n}}\left(i^{\mbox{\emph{{tr}}}[r\alpha x]}\right)_{x\in\mathcal{T}_r}
\end{align}
and hence $r\star \vec{v}_\alpha\in M$, for all $r\in GR(4,1)$.  Hence $M$ is a submodule. 

2.  Part 1. shows that $M$ is a module.  To show $M$ is free we need that for every $\vec{v}$ such that $2\vec{v}=0$, there exists $\vec{u}$ such that $2\vec{u}=\vec{v}$.  Thus we require that if $\alpha$ is such that 
\begin{equation}
2\textrm{tr}(\alpha x)=\textrm{tr}(2\alpha x)=0
\end{equation}
for all $x\in \mathcal{T}_n$, then there exists $\beta\in M$ such that $\alpha=2\beta$.  Reverting to the p-adic notation, let $\alpha=a+2b$ and $\beta=c+2d$, then $2\alpha=0+2a$ and $2\beta=0+2c$.  Hence we need to show that if $\textrm{tr}(2ax)=0$ for all $x\in\mathcal{T}_n$, then $a=0$.

Using Theorem \ref{thm:trace}.\ref{part:0}, we see that this is equivalent to showing that for all $x\in\mathcal{T}_n$, there exists $\gamma=(e+2f)\in GR(4,n)$ such that \begin{align}
2ax & = e+2f-\phi(e+2f)\\
2ax & =e+2f-e^2-2f^2
\end{align}
where $a,x,e,f\in\mathcal{T}_n$. This simplifies to 
\begin{equation}\label{eqn:ff}
ax  = f-f^2
\end{equation}
If $a=0$, then we have solved our problem. Assume $a\neq 0$, then there exists $x\in \mathcal{T}_n$ such that $ax=1$.  Thus we require a solution to 
\begin{equation}\label{eqn:f}
 0 = f^2-f+1
\end{equation}
This is a monic irreducible polynomial of degree $2$, and hence has possible solution only in $GR(4,2)$.  Let $h(f)=f^2-f+1$, then $GR(4,2)=\mathbb{Z}_4[f]\slash (h(f))$, and hence $\mathcal{T}_2=\{0,1,\xi,\xi+3\}$ where $\xi$ is a root of $h(f)$. From equation (\ref{eqn:ff})
\begin{align}
\xi-\xi^2= & \xi-\xi-3=1,\\
\xi^2-\xi^4= & \xi^2-\xi= 3.
\end{align}
Hence if $ax\in\{\xi,\xi+3\}$ then  equation (\ref{eqn:ff}) has no solution.  We require that equation (\ref{eqn:ff}) holds for fixed $a$ and all $x\in\mathcal{T}_n$, hence we require that $a=0$, which shows that    $M$ is a free submodule. 
%Checking the specific case of $n=2$, using equation (\ref{eqn:f}), $\mathcal{T}_2=\{0,1,\xi,\xi+3\}$.  
%\begin{align}
%\textrm{tr}(2a0) & = 0\\
%\textrm{tr}(2a) & = 2(a+a^2)=2 \quad \forall a\in \mathcal{T}_n\setminus 0\\
%\textrm{tr}(2a\xi)& =2(a\xi+a^2\xi^2)=2\\
%\textrm{tr}(2a\xi^2)& =2(a\xi^2+a^2\xi)=2
%\end{align}
And thus by construction forms a subspace of $PHG(2^n-1,GR(4,1))$.  The counting results of Lemma \ref{lem:evencount} show the size of the subspace.
\end{IEEEproof}
Note that $GR(p^s,1)\cong \mathbb{Z}_{p^s}$, and as with Theorem \ref{thm:odd}, $M=M'$ and $N=N'$,  thus the conditions of Proposal \ref{prop} are satisfied.
\section{Conclusion}
We have shown that several sets of MUBs display the algebraic structure of a module  and the geometric structure of a subspace of a projective Hjelmslev geometry.  There are also counting results to show that this geometric structure may be true in general.
Of particular note is that these structures may not arise from the sets of vectors which define the MUBs, but from the sets of vectors derived from component wise multiplication.

We have not covered all possible constructions of MUBs, but have shown sufficient evidence that this is a structure worthy of more thorough investigation.

\end{document}